\documentclass[11pt]{article}
\usepackage{amsmath}
\usepackage{amsmath,amssymb,amsfonts,amsthm,fancyhdr}
\usepackage{epsfig,graphicx,picins,picinpar,subfigure}
\usepackage{pstricks}
\usepackage{fancyvrb}
\usepackage[numbers,sort&compress]{natbib}

\begin{document}

\title{\textbf{Goldbach Conjecture and the least prime number in an arithmetic
progression}}
\author{Shaohua Zhang}
\date{{\small  \emph{School of Mathematics, Shandong
University, Jinan, Shandong, 250100, China}}}

\footnotetext [1]{This work was partially supported by the National
Basic Research Program (973) of China (No. 2007CB807902) and the
Natural Science Foundation of Shandong Province (No. Y2008G23).}
\footnotetext [2]{E-mail address: shaohuazhang@mail.sdu.edu.cn}
\maketitle

\vspace{3mm}\textbf{Abstract}

\vspace{3mm}In this Note, we try to study the relations between the
Goldbach Conjecture and the least prime number in an arithmetic
progression. We give a new weakened form of the Goldbach Conjecture.
We prove that this weakened form and a weakened form of the Chowla
Hypothesis imply that every sufficiently large even integer may be
written as the sum of two distinct primes.

\vspace{3mm}\textbf{R\'{e}sum\'{e}}

\vspace{3mm}\textbf{La conjecture de Goldbach et le plus petit
nombre premier dans une progression arithm\'{e}tique}

\vspace{3mm} Dans ce document, nous essayons d'\'{e}tudier les
relations entre la conjecture de Goldbach et le plus petit nombre
premier dans une progression arithm\'{e}tique. Nous donnons une
nouvelle forme faible de la conjecture de Goldbach. Nous prouvons
que cette forme affaiblie et une forme affaiblie de l'hypoth\`{e}se
de Chowla impliquent que tout entier pair suffisamment grand peut
\^{e}tre \'{e}crit comme une somme de deux nombres premiers
distincts.

\vspace{3mm}\textbf{Keywords:} Goldbach Conjecture, least prime
number, arithmetic progression,  Chowla Hypothesis, Generalized
Riemann Hypothesis

\vspace{3mm}\textbf{2000 MR  Subject Classification:} 11A41, 11A99,
11B25, 11P32


\section{Introduction}
\setcounter{section}{1}\setcounter{equation}{0}
Goldbach's famous conjecture states that every even integer $2n\geq
4$ is the sum of two primes. Since it is trivial that for infinitely
many even integers: $2p=p+p$ (for every prime $p$), we give a
slightly different form of this conjecture: every even integer
$2n\geq 8$  is the sum of two distinct primes. Thus, one can state
Conjecture 1 below, which is also called a weakened form of
Goldbach's Conjecture or the necessary condition of Goldbach's
Conjecture.

\vspace{3mm}\noindent{\bf  Conjecture 1.~~}%
For integer $n>5$, there exists a natural number $r$  such that
$2n-p_r$ is coprime to each of
$2n-p_1,...,2n-p_{r-1},2n-p_{r+1},...,2n-p_k$, where
$p_1,...,p_{r-1},p_r, p_{r+1},...,p_k$ are all old primes smaller
than $n$, $p_r$ satisfies $(p_r,n)=1$ and $1\leq r\leq k=\pi
(n-1)-1$.

\vspace{3mm}Let  $k,l$ denote positive integers with $(k,l)=1$ and
$1\leq l\leq k-1$. Denote by $p(k,l)$ the least prime $p\equiv
l(\mod k)$. Let $p(k)$ be the maximum value of $p(k,l)$ for all $l$
with $(k,l)=1$ and $1\leq l\leq k-1$. In 1992, Heath-Brown [2]
proved $p(k)\ll k^{5.5}$. This is the  best known result on  $p(k)$.
Recently, Heath-Brown told the author that  Xylouris
(http://arxiv.org/abs/0906.2749) has improved his result to $p(k)\ll
k^{5.2}$.   Chowla [1] has observed that $p(k)\ll k^{2+\epsilon}$
for every $\epsilon>0$ assuming the Generalized Riemann Hypothesis.
He further conjectured $p(k)\ll k^{1+\epsilon}$ for every
$\epsilon>0$. Based on the conjecture of Chowla, one might state the
following, Conjecture 2:

\vspace{3mm}\noindent{\bf  Conjecture 2.~~}%
For every sufficiently large positive integer $k$, namely when
$k>c_1$,  $p(k)<k^{1.5}$, where  $ c_1$ is a positive constant.

\vspace{3mm} The object of this Note is to study the relations
between the Goldbach Conjecture and the least prime number in an
arithmetic progression. We obtained the following Theorem 1 which
gives a sufficient condition for the Goldbach Conjecture. As we
know, even under Riemann hypothesis or if the generalized Riemann
hypothesis holds, nobody has proved up until now that the Goldbach
Conjecture is true. Therefore, needless-to-say, that refining the
results of Heath-Brown and Xylouris, and proving Conjecture 1,
should be given much attention.

\vspace{3mm}\noindent{\bf  Theorem  1.~~}%
\emph{If Conjecture 1 and Conjecture 2 hold, then every sufficiently
large even integer may be written as the sum of two distinct
primes.}

\section{The proof of Theorem 1}
\vspace{3mm}\noindent{\bf  Proof.~~}%
One can prove that for every prime $p\geq 48673$,  and any integer
$a$ with $1\leq a< p^{1.5}$, there is a prime $q$ coprime to $a$ and
such that $4q^3<p$.

By the prime number theorem in an arithmetic progression, it is easy
to prove that for any prime $p$ with $p\leq\max\{c_1,48673\}$,
($c_1$ is the positive constant in Conjecture 2), there exists a
positive constant $c_2>6$ such that for every positive integer
$n>c_2$, when $(p,n)=1$, there exist two distinct odd primes $p_1$
and $p_2$ satisfying $2n\equiv p_1\equiv p_2(\mod p)$ and
$p_1,p_2\in Z_n^*=\{x|1\leq x\leq n, (x,n)=1\}$.

Let $n$ be an integer $>c_2$. Since we assume Conjecture 1, there
exists  $r>1$ such that $(p_r,n)=1$ and $2n-p_r$ is coprime to every
$2n-p$ when $p$ ranges through the odd primes $\leq n$ and different
from $p_r$. We will show that $2n-p_r$ is prime. If this is the
case, then Theorem 1 is proved, so let us suppose we can write
$2n-p_r =pm$, where $p$ is the least prime factor of $2n-p_r$. Thus,
$2n>p^2$.

We have $p>\max\{c_1,48673\}$. Indeed, if $p$ is smaller, we can
find two odd primes say $q_1$ and $q_2$, not more than $n$ and prime
to $2n$, such that $2n\equiv q_1\equiv q_2(\mod p)$. At most one of
them, say $q_1$, can be equal to $p_r$. This means that $2n-p_r$ is
not coprime to $2n-q_2$, contrarily to our hypothesis on $p_r$.

Note that $p_r\neq p$ since $(p_r,n)=1$. If $p_r<p$, then
$p+p_r<p^{1.5}$ and there is a prime $q$ coprime to $p+p_r$ and such
that $4q^3<p$. Since we suppose that Conjecture 2 holds, hence there
is a prime $x$ such that $x\equiv p+p_r(\mod pq)$ and
$x<(pq)^{1.5}<\frac{p^2}{2}<n$. Clearly, $p_r\neq x$. But
$p|(2n-p_r,2n-x)$. It is a contradiction by our assumption on $p_r$.

Hence $p_r>p$. We write $p_r=pl+v$ with $1\leq v<p$. If  $l\geq
\sqrt{p}$, there is a prime $y$ such that $y\equiv v(\mod p)$ and
$y<p^{1.5}<p_r$ (since we suppose Conjecture 2). But we have also
$p|(2n-p_r,2n-y)$, it is contrary to our assumption on $p_r$ again.
So we have $l< \sqrt{p}$, $lv<p^{1.5}$ and there is a prime $q$
coprime to $lv$ and such that $4q^3<p$. Note that there is a prime
$z$ such that $z\equiv v(\mod pq)$ and
$z<(pq)^{1.5}<\frac{p^2}{2}<n$ (since we suppose that Conjecture 2
holds). Obviously, we have $z\neq p_r$ since $(q,l)=1$. But
$p|(2n-p_r,2n-z)$. The contradiction implies that $2n-p_r$ is a
prime number. This completes the proof of Theorem 1.

\section{Acknowledgements}
I am very greatfull to the referees and Professor Heath-Brown for
their comments improving the presentation of the Note, and also to
my supervisor Professor Xiaoyun Wang for her suggestions. Thanks
also go to Mingqiang Wang, Huaning Liu and Hongbo Yu for their help.
I thank my father and my wife for their encouragement and also thank
the key lab of cryptography technology and information security in
Shandong University and the Institute for Advanced Study in Tsinghua
University, for providing me with excellent conditions.

\clearpage
\end{document}